\documentclass{elsart}

\usepackage{graphicx}
\usepackage{amssymb}
\usepackage[T1]{fontenc}
\usepackage{latexsym}
\usepackage{xypic}
\usepackage{eufrak}
\usepackage{euscript}
\usepackage{amsfonts,amsmath,mathrsfs}
\usepackage{verbatim}
\usepackage{fancyhdr}
\usepackage[english]{babel}

\newtheorem{theo}{Theorem}
\newtheorem{defi}[theo]{Definition}

\newtheorem{pro}[thm]{Proposition}

\def\cal{\mathcal}
\renewcommand{\c}{\widehat}
\renewcommand{\b}{\overline}
\renewcommand{\d}{\displaystyle}
\newcommand{\ep}{\varepsilon}
\renewcommand{\geq}{\geqslant}
\renewcommand{\leq}{\leqslant}
\def\scr{\mathscr}

\newcommand{\N}{\mathbb{N}}

\newcommand{\R}{\mathbb{R}}
\newcommand{\Q}{\mathbb{Q}}
\newcommand{\C}{\mathbb{C}}

\newcommand{\bP}{\mathbb{P}}

\def\cal{\mathcal}

\def\oT{{[0{,}T]}}

\def\1{{\mathbf{1}}}

\def\div{{\rm{div}}}

\def\sg{\sigma}

\begin{document}

\begin{frontmatter}



\title{Dynamical properties and characterization of gradient drift diffusions}

\author[label1]{S\'ebastien Darses}
\author[label1]{and Ivan Nourdin}
 \address[label1]{Université Pierre et Marie Curie Paris VI\\
Laboratoire de Probabilit\'es et Mod\`eles Al\'eatoires\\
Bo\^ite courrier 188, 4 place Jussieu,
75252 Paris Cedex 05, France\\
{\tt \{sedarses,nourdin\}@ccr.jussieu.fr}
}

\begin{abstract}
We study the dynamical properties of the Brownian diffusions having
$\sigma\,{\rm Id}$ as diffusion coefficient matrix and $b=\nabla U$
as drift vector. We characterize this class through the equality
$D^2_+=D^2_-$, where $D_{+}$ (resp. $D_-$) denotes the forward
(resp. backward) stochastic derivative of Nelson's type. Our proof
is based on a remarkable identity for $D_+^2-D_-^2$ and on the use
of the martingale problem. We also give a new formulation of a
famous theorem of Kolmogorov concerning reversible diffusions. We
finally relate our characterization to some questions about the
complex stochastic embedding of the Newton equation which initially
motivated of this work.
\end{abstract}

\end{frontmatter}

% \appendix

% \section{}
% \label{}

\newpage

\section{Introduction}
\label{intro}

%In this paper, we are interested in the dynamical properties of
%gradient drift diffusions with a constant diffusion coefficient.
%In other words, we will study the
%stochastic derivatives of these diffusions. Recall that, f
For a general process
$Z=(Z_t)_{t\in[0,T]}$ defined on a probability space
$(\Omega,\mathscr{F},\bP)$, we have introduced in \cite{dn} the notion of {\it stochastic
derivative} for $Z$ at $t$ with respect to a {\it differentiating} sub-$\sg$-field ${\scr A}^t$ of $\scr{F}$
(resp. {\it forward differentiating}, {\it backward differentiating}). More precisely, it
means that ${\scr A}^t$ is such that the quantity
$${\rm E}\left[\frac{Z_{t+h}-Z_t}{h} |{\scr A}^t\right]$$
converges in probability (or for another topology) when $h\rightarrow 0$
(resp. $h\downarrow 0$, $h\uparrow 0$); the limit being called
the stochastic derivatives of $Z$ at $t$ w.r.t. $\scr{A}^t$. When we consider Brownian
diffusions of the form
\begin{equation}\label{X}
X_t=X_0+\int_0^tb(s,X_s)ds+\int_0^t\sg(s,X_s)dW_s,\quad t\in[0,T],
\end{equation}
then, under suitable conditions, the $\sg$-field $\scr{T}^X_t$ generated by $X_t$ is
both a forward and backward differentiating $\sg$-field for $X$ at $t$. The
associated derivatives are called {\it Nelson derivatives}, due to the
Markov property of the diffusion and of its time reversal which allow
to take the conditional expectation both with respect to the
past $\scr P^X_t$ and the future $\scr F^X_t$ of the diffusion.
For simplicity, we note them respectively $D_+$ and $D_-$
in the sequel.
Notice that these derivatives are relevant and natural quantities for Brownian
diffusions: they are indeed respectively equals to the forward and the backward (up to sign) drift of $X$.
Moreover, they exist under rather mild conditions, see
{\it e.g.} \cite{mns, pardoux}.

We shall see that Nelson derivatives turn out to have remarkable
properties when we work with diffusions of the type
\begin{equation}\label{constant}
X_t=X_0+\int_0^t b(s,X_s)ds+\sigma\,W_t,\quad t\in[0,T].
\end{equation}
Here, $\sigma\in\R$ is assumed to be constant. For instance, we
shall show that the equalities $D_+X_t=-D_-X_t$, $t\in(0,T)$,
characterizes the class of stationary diffusions of the type
(\ref{constant}) having moreover an homogeneous gradient drift (see
Proposition \ref{firstorder}). This statement is in fact quite easy
to obtain. A more difficult one, which is the main result of this
paper, states that a Brownian diffusion of the type (\ref{constant})
is a gradient diffusion - that is, its drift coefficient writes
$b=\nabla_x U$ for a certain $U$ - if and only if
$D^2_+X_t=D^2_-X_t$ for any $t\in(0,T)$, see Theorem \ref{main} for
a precise statement. Let us notice that this result was conjectured
at the end of the note \cite{cd}. Our proof is based on the
discovery of a remarkable identity (Lemma \ref{div}): we can write
the quantity $p_t(X_t)(D^2_+X_t-D^2_-X_t)$ as the divergence of a
certain vector field, where $p_t$ denotes the density of the law of
$X_t$. Combined with the expression of the adjoint of the
infinitesimal generator, we can then conclude using probabilistic
arguments, especially the martingale problem. Let us moreover stress
on the fact that we were able to solve our problem with
probabilistic tools, whereas its analytic transcription with the
help of partial differential equations seemed more difficult to
treat.

The paper is organized as follows. In section 2, we introduce some
notations and we give the useful expressions of the Nelson
derivatives under the conditions given by Millet, Nualart and Sanz
in \cite{mns}. In section 3, we study the above mentioned
characterizations and we prove our main result. In section 4, we
make some remarks on the questions related to the complex stochastic
embedding of the Newton equation, which have motivated this work.

\section{Preliminaries on stochastic derivatives}
\label{}

\subsection{Notations}\label{s21}

Let $T>0$ and $d\in\N^*$. The space $\R^d$ is endowed with its
canonical scalar product $\langle\cdot,\cdot\rangle$. Let $|\cdot|$
be the induced norm.

If $f:\oT\times \R^d\to \R$ is a smooth function, we set
$\partial_jf=\frac{\partial f}{\partial x_j}$. We denote by $\nabla
f=(\partial_i f)_i$ the gradient of $f$ and by $\Delta f=\sum_j
\partial_j^2f$ its Laplacian. For a smooth map
$\Phi:\oT\times \R^d\to\R^d$, we denote by $\Phi^j$ its
$j^{th}$-component, by $(\partial_x \Phi)$ its differential which we
represent into the canonical basis of $\R^d$: $(\partial_x
\Phi)=(\partial_j\Phi^i)_{i,j}$, and by $\div \Phi=\sum_j
\partial_j\Phi^j$ its divergence. By convention, we denote by
$\Delta \Phi$ the vector $(\Delta \Phi^j)_j$. The image of a vector
$u\in\R^d$ under a linear map $M$ is simply denoted by $Mu$, for
instance $(\partial_x \phi)u$. The map
$a:\oT\times\R^d\to\R^d\otimes\R^d$ is viewed as $d\times d$
matrices whose columns are denoted by $a_k$. Finally, we denote by
$\div\,a$ the vector $(\div \,a_k)_k$.

Let $(\Omega,\mathscr{A},\bP)$ be a probability space on which is
defined a $d$-dimensional Brownian motion $W$. For a process $Z$
defined on $(\Omega,\mathscr{A},\bP)$, we set $\scr P_t^Z$ the $\sigma$-field generated by $Z_s$ for $0\leq s\leq t$ and $\scr F_t^Z$ the
$\sigma$-field generated by $Z_s$ for $t\leq s\leq T$. Consider the
$d$-dimensional diffusion process $X=(X_t)_{t\in[0,T]}$ solution of
the stochastic differential equation (\ref{X}) where $X_0\in{\rm
L}^2(\Omega)$ is a random vector independent of $W$, and the functions
$\sigma:\oT\times\R^d\to\R^d\otimes\R^d$ and
$b:\oT\times\R^d\to\R^d$ are Lipschitz with linear growth.
More precisely, we assume that $\sigma$ and $b$ satisfy the two
following conditions: there exists a constant $K>0$ such that, for
all $x,y\in\R^d$, we have
    $$\sup_{t\in[0,T]} \big[\left|b(t,x)-b(t,y)\right|+\left|\sigma(t,x)-\sigma(t,y)\right|\big]
\leq K\left|x-y\right|$$ and
$$\sup_{t\in[0,T]}\big[\left|b(t,x)\right|+\left|\sigma(t,x)\right|\big]\leq
    K(1+\left|x\right|).$$

%We denote by $\mathscr{C}^k_b$ the set of all bounded functions $f:\oT\times
%\R^d\to\R^d$ such that all their derivatives up to order $k$ are
%also bounded.
We moreover assume that $b$ is differentiable w.r.t. $x$ and we set
$G=(\partial_x b)-(\partial_x b)^*$, \textit{i.e.}
$G_i^j=\partial_ib^j-\partial_jb^i$. Finally, we set
$a=\sigma\sigma^*$, {\it i.e.} $a_i^j=\sum_{k}\sigma_i^k\sigma_j^k$.

%and an independent r.v. $X_0$.
In the sequel, we will work under the following assumption:
\begin{itemize}
\item[(H)] For any $t\in (0,T)$, the law of $X_t$ admits a positive
density $p_t:\R^d\to (0,+\infty)$ and we have, for any $t_0\in
(0,T)$:
\begin{equation}\label{H_1}
\max_{j=1,\ldots,n}\int_{t_0}^T
\int_{\R^d} \left|\div(a_j(t,x)p_t(x))\right|dxdt <
    +\infty.
\end{equation}
The functions
\begin{equation}\label{H_2}
\frac{\div(a_j(t,\cdot)p_t(\cdot))}{p_t(\cdot)}
\end{equation}
are Lipschitz.
\end{itemize}

The condition (\ref{H_1}) will ensure us that the time reversed
process $\b X_t=X_{T-t}$ is again a diffusion process (see
\cite{mns}, Theorem 2.3). Let us moreover notice that our condition (\ref{H_2})
is weaker than that which is imposed in Proposition 4.1 of \cite{T}. 
Finally, let us remark that the positivity
assumption made on $p_t$ is quite weak when $X$ is of the type
(\ref{constant}): it is for instance automatically verified when we
can apply Girsanov theorem in (\ref{constant}), that is when the
Novikov condition is verified.

\subsection{Stochastic derivatives of Nelson's type}

In \cite{dn}, we have introduced the notion of {\it differentiating}
$\sg$-field:

\begin{defi}\label{defi1}
    Set $t\in(0,T)$ and let $Z$ be a process. We say that $\scr A^t$ (resp. $\scr B^t$) is
    a {\rm forward differentiating }$\sigma$-{\rm field} (resp. {\rm backward differentiating}
    $\sigma$-{\rm field}) for $Z$ at $t$ if $E[\frac{Z_{t+h}-Z_t}{h}|\scr A^t]$
(resp. $E[\frac{Z_{t}-Z_{t-h}}{h}|\scr B^t]$) converges in
probability when $h\downarrow 0$. In these cases, we define the
so-called forward and backward derivatives
\begin{eqnarray}
D^{\scr A^t}_+ Z_t & = & \lim_{h\downarrow 0}E\left[
\frac{Z_{t+h}-Z_t}{h}\,|\scr A^t\right],\label{d_forward}\\
D^{\scr B^t}_- Z_t & = & \lim_{h\downarrow 0}E\left[
\frac{Z_{t}-Z_{t-h}}{h}\,|\scr B^t\right] \label{d_backward}.
\end{eqnarray}
\end{defi}

For Brownian diffusions $X$ of the form (\ref{X}), the
present turns out to be a forward and backward differentiating
$\sg$-field. Precisely, the $\sg$-field ${\scr T}^X_t$ generated by $X_t$ is
both forward and backward differentiating for $X$ at $t$. Equivalently,
due to the Markov property of $X$ (resp. of its time reversal $\b X$),
$\scr P^X_t$ (resp. $\scr F^X_t$) is forward (resp. backward)
differentiating for $X$ at $t$. For this reason, we call the
derivatives defined by (\ref{d_forward}) and (\ref{d_backward})
stochastic derivatives of Nelson's type. Indeed, in \cite{n} Nelson
introduced the processes which have stochastic derivatives in
$L^2(\Omega)$ with respect to a fixed filtration $(\scr
P_t)$ and a fixed decreasing filtration $(\scr F_t)$.

Henceforth, we work with the stochastic derivatives of Nelson's type
for Brownian diffusions and so we simply write $D_{\pm}X$ instead of
$D_{\pm}^{{\scr T}^X_t}X_t$. Now, we can relate the stochastic
derivatives of Nelson's type to the time reversal theory:

\begin{pro}\label{derivgood}
Let $X$ be given by (\ref{X}) and satisfying assumption (H). Then $X$
is a Markov diffusion w.r.t. the increasing filtration
$(\scr{P}^X_t)$ and the decreasing filtration $(\scr{F}^X_t)$.
Moreover for almost all $t\in(0,T)$, $\scr{T}_t^X$ is a forward and
backward differentiating $\sg$-field for $X$ at $t$ and
\begin{eqnarray}
D_+X_t & = & b(t,X_t)\\
D_-X_t & = & b(t,X_t)-\frac{\div(a(t,X_t)p_t(X_t))}{p_t(X_t)}.
\end{eqnarray}
\end{pro}

\begin{pf*}{Proof.}
The proof essentially uses  Theorem 2.3 of Millet-Nualart-Sanz
\cite{mns}, and is divided in two steps:

1) $X$ is a Markov diffusion w.r.t. the increasing filtration
$(\scr{P}^X_t)$, so:
$$E\left[\frac{X_{t+h}-X_t}{h}\left|\scr{P}^X_t\right.\right]=E\left[\frac{1}{h}\int_t^{t+h}b(s,X_s)ds\left|\scr{P}^X_t\right.\right],
$$ and
\begin{eqnarray*}
  E\left|E\left[\left.\frac{X_{t+h}-X_t}{h}\right|\scr{P}^X_t\right]-b(t,X_t)\right| & \leq &
  \frac{1}{h} E\int_t^{t+h}\left|b(s,X_s)-b(t,X_t)\right|ds \\
    &=& \frac{1}{h}
    \int_t^{t+h}E\left|b(s,X_s)-b(t,X_t)\right|ds.
\end{eqnarray*}
Using the fact that $b$ is Lipschitz and that $t\mapsto E|X_t|$ is
locally integrable (see, {\it e.g.}, Theorem 2.9 in \cite{KS}), we
can conclude by the differentiation Lebesgue theorem that for almost
all $t\in(0,T)$:
$$\d \frac{1}{h}\int_t^{t+h}E\left|b(s,X_s)-b(t,X_t)\right|ds \rightarrow 0\mbox{ a.s.,}\quad\mbox{as }h\to 0.$$

Therefore $D_+X_t$ exists and is equal to $b(t,X_t)$.

2) Thanks to assumption (H), we can apply Theorem $2.3$ in
\cite{mns}. Hence $\overline{X}_t=X_{T-t}$ is a diffusion process
w.r.t. the increasing filtration $(\scr{F}_{T-t})$ and whose
generator reads $$\overline{L}_tf=\overline{b}^i\partial_i
f+\frac{1}{2}\overline{a}^{ij}\partial_{ij}f$$ with
$\overline{a}^{ij}(T-t,x)=a^{ij}(t,x)$ and
$$
\overline{b}^{i}(T-t,x)=-b^i(t,x)+\frac{\div(a_i(t,x)p_t(x))}{p_t(x)}.$$
We have :
\begin{eqnarray}
    E\left[\left.\frac{X_t-X_{t-h}}{h}\right|\scr{F}^X_t \right] & = & E\left[\left.\frac{\overline{X}_{T-t}-
    \overline{X}_{T-t+h}}{h}\right|\scr{F}^X_{T-t}\right] \nonumber\\
    & = & -E\left[\left.\frac{1}{h}\int_{T-t}^{T-t+h}\overline{b}(s,\overline{X}_s)ds\right|\scr{F}^X_{T-t}\right].
\end{eqnarray}
Assumption (H) implies that $$t\mapsto
E\left|\frac{\div(a_i(t,X_t)p_t(X_t))}{p_t(X_t)}\right|$$ is locally
integrable. Then, using the same calculations and arguments as
above, we obtain that $D_-X_t$ exists and is equal to
$-\overline{b}(T-t,\overline{X}_{T-t})$. \qed
\end{pf*}

\begin{cor}\label{cor_sigma_constant}
If $X$ given by (\ref{constant}) verifies assumption (H), we have for
almost all $t\in (0,T)$:
$$
  D_+X_t = b(t,X_t) \quad\mbox{ and }\quad
  D_- X_t = b(t,X_t)-\sg^2 \frac{\nabla p_t}{p_t}(X_t) .
$$
\end{cor}

\begin{rem}
{\rm The appearance of the density $p_t$ in the formula giving
$D_-X_t$ may seem surprising at first sight. As clear from the proof
of Theorem \ref{derivgood}, the reason for such a formula stems from
 the Brownian theory of time reversal. The same term was obtained
by F\"{o}llmer \cite{follmer} for Brownian semimartingales of the
form $\int_0^tb_sds+W_t$ with $E\int_0^Tb_s^2ds<\infty$, by relating
the backward Nelson derivative and the time reversed drift. Based on
the same strategy, Millet, Nualart and Sanz \cite{mns} extended the
result to diffusions satisfying (H) using Malliavin calculus.
Finally, this additional term can also be viewed as the result of a
"{\it grossissement de filtration}" (see Pardoux \cite{pardoux}).
Roughly speaking, when we consider a diffusion
$X_t=\int_0^tb(s,X_s)ds+\int_0^t\sigma(s,X_s)dW_s$ and $\scr G_t$
the $\sigma$-field generated by $W_u-W_r$ for $T-t\leq u<r\leq T$,
then $\b W_t-\b W_0$ is a $\scr G_t$-Brownian motion and the
question sums up to writing the Doob-Meyer decomposition of $\b
W_t-\b W_0$ in the enlarged filtration $\scr H_t=\scr G_t\vee \b
X_t$. In particular, knowing this answer gives the decomposition of
$\b X$ with respect to its natural filtration.}
\end{rem}

Finally, we will also need the following composition formula, stated
by Nelson \cite{n} and that we prove for the diffusions we consider.

\begin{pro}\label{nelson-composition}
Let $f\in C^{1,2}([0,T]\times\R^d)$ with bounded second order
derivatives and let $X$ be a diffusion of the form (\ref{constant})
satisfying (H). Then, for almost all $t\in(0,T)$:
\begin{equation}\label{}
    D_{\pm}f(t,X_t)=\left(\partial_t f+ (\partial_xf) D_{\pm}X_t
    \pm\frac{\sg^2}{2}\Delta f\right)(t,X_t).
\end{equation}

\end{pro}
\begin{pf*}{Proof.} Let $h>0$.\\
1) The forward case. The Taylor formula yields:
\begin{eqnarray}
  f(t+h,X_{t+h})-f(t,X_t) &\ =\ & \partial_tf(t,X_t)h+ \partial_x f(t,X_t)(X_{t+h}-X_t) \label{dev}\\
   & & +\frac{1}{2}\sum_{i,j=1}^n(X^i_{t+h}-X^i_t)(X^j_{t+h}-X^j_t)\partial^2_{ij}f(t,X_t) +
R(t,h) \nonumber
\end{eqnarray}
where the remainder $R(t,h)$ is given by
\begin{eqnarray*}
R(t,h)&=&\frac{1}{2}\sum_{i,j=1}^n(X^i_{t+h}-X^i_t)(X^j_{t+h}-X^j_t)
\big(\partial^2_{ij}f(u_{t,h})
- \partial^2_{ij}f(t,X_t)\big)\\
& & +h\,\sum_{j=1}^n(X^j_{t+h}-X^j_t)
\partial_t\partial_{j}f(u_{t,h})
\end{eqnarray*}
with $u_{t,h}=(t+\theta h,(1-\theta)X_t+\theta X_{t+h})$ and
$\theta\in (0,1)$ depending on $t$ and $h$.

We first treat the third term of the r.h.s of (\ref{dev}). For
instance for the term $\frac{1}{h}E[(X^i_{t+h}-X^i_t)^2|X_t]$:
\begin{equation}\label{(X-X)}
(X^i_{t+h}-X^i_t)^2=\left(\int_t^{t+h}b(s,X_s)ds\right)^2+\sg^2(W^i_{t+h}-W^i_t)^2+2\sg(W^i_{t+h}-W^i_t)\int_t^{t+h}b(s,X_s)ds.
\end{equation}
We have by Schwarz inequality:
$$\left(\int_t^{t+h}b(s,X_s)ds\right)^2\leq h \int_t^{t+h}b^2(s,X_s)ds.$$
Thus $$\frac{1}{h}E\left(\int_t^{t+h}b(s,X_s)ds\right)^2\leq
\int_t^{t+h}E[b^2(s,X_s)]ds \longrightarrow 0,$$
since $t\rightarrow E|X_t|^2$ is locally integrable
(see, {\it e.g.}, Theorem 2.9 in \cite{KS}). Again by Schwarz
inequality, we deduce that
$h^{-1}\big(W^i_{t+h}-W^i_t\big)\int_t^{t+h}b(s,X_s)ds$ tends to $0$ in
$L^1(\Omega)$. Moreover:
$$\frac{1}{h}E[(W^i_{t+h}-W^i_t)^2|X_t]=\frac{1}{h}E[(W^i_{t+h}-W^i_t)^2]=1.$$

We now treat the remainder of (\ref{dev}). The fact that $\partial^2
f$ is bounded allows to show as above that
$h^{-1}\left(\int_t^{t+h}b(s,X_s)ds\right)^2(\partial^2_{ij}f(u_{t,h})
- \partial^2_{ij}f(t,X_t))$ and
$$\frac{W^i_{t+h}-W^i_t}{h}\int_t^{t+h}b(s,X_s)ds(\partial^2_{ij}f(u_{t,h})-
\partial^2_{ij}f(t,X_t))$$ converges to $0$ in $L^1(\Omega)$.
Moreover
$$
\begin{array}{lll}
E\left[\frac{(W^i_{t+h}-W^i_t)^2}{h}(\partial^2_{ij}f(u_{t,h}) -
\partial^2_{ij}f(t,X_t))\right]\\
\quad\quad\quad\quad\quad\quad\quad\quad
\le
\frac{\sqrt{E|W^i_{t+h}-W^i_t|^4}}{h}\sqrt{E|\partial^2_{ij}f(u_{t,h})
- \partial^2_{ij}f(t,X_t)|^2}\\
\quad\quad\quad\quad\quad\quad\quad\quad
\le
C \sqrt{E|\partial^2_{ij}f(u_{t,h})
- \partial^2_{ij}f(t,X_t)|^2}.
\end{array}
$$
Since $\partial^2 f$
is bounded and $u_{t,h}$ tends to $(t,X_t)$ as $h\rightarrow 0$, we can apply the bounded
convergence theorem and conclude.

2) The backward case. We calculate the Taylor expansion of
$-(f(t-h,X_{t-h})-f(t,X_t))$ and we write $(X^i_{t-h}-X^i_t)^2=(\b
X^i_{T-t+h}-\b X^i_{T-t})^2$. We then write the decomposition
(\ref{(X-X)}) for $\b X$ with its time reversed drift $\b b$ and its
time reversed driving Brownian motion $\c W$. So the computations
are identical to those of the first point. \qed
\end{pf*}

\section{Dynamical study of gradient diffusions}
\label{}

\subsection{First order derivatives}

In this section, we only consider Brownian diffusions of type
(\ref{constant}) with a homogeneous drift. More precisely, we work
with $X$ verifying
\begin{equation}\label{constant+hom}
X_t=X_0+\int_0^t b(X_s)ds+\sigma\,W_t,\quad t\in[0,T].
\end{equation}
We can then characterize the sub-class of stationary diffusions
having a gradient drift vector, by means of first order Nelson
derivatives. Such diffusions were already considered by many
authors. A result of Kolmogorov \cite{kolmo} states that $b$ is a
gradient if and only if the law of $X$ given by (\ref{constant+hom})
is reversible, {\it i.e.} $(X_t)_{t\in[0,T]}$ and
$(X_{T-t})_{t\in[0,T]}$ have the same law. In what follows, we show
that another characterization of this last fact can be made with the
help of Nelson derivatives. For instance, knowing that $b$ is a
gradient allows to easily construct an invariant law for $X$. More
precisely, when $b=\nabla U$ with $U:\R^d\to\R$ regular enough and
with sufficiently fast decrease at infinity, the probability law
$\mu$ defined by
$$
d\mu=c^{-1}{\rm e}^{\frac{2U(x)}{\sigma^2}}dx\quad\mbox{with}\quad
c=\int_{\R^d}{\rm e}^{\frac{2U(x)}{\sigma^2}}dx<\infty
$$ is
invariant for $X$.

We can easily prove  the following:
\begin{pro}\label{firstorder}
Let $X$ be the Brownian diffusion defined by (\ref{constant+hom}). We moreover assume
that $X$ verifies assumption (H).
\begin{enumerate}
\item
If $D_+X_t=-D_-X_t$ for any $t\in(0,T)$ then $b=\nabla U$ with
$U:\R^d\to\R$ given by $U=\frac{\sigma^2}{2}\log p_t$. In
particular, $X$ is a stationary diffusion with initial law $\mu$
given by $d\mu={\rm e}^{\frac{2U(x)}{\sigma^2}}dx$.
\item
Conversely, if $b=\nabla U$ with $U:\R^d\to\R$ such that
$c:=\int_{\R^d}e^{\frac{2U(x)}{\sigma^2}}dx<\infty$ and if the law
of $X_0$ is $d\mu=c^{-1}{\rm e}^{\frac{2U(x)}{\sigma^2}}dx$, then
the probability law $\mu$ is invariant for $X$ and, for any
$t\in(0,T)$, we have $D_+X_t=-D_-X_t$.
\end{enumerate}
\end{pro}
\begin{pf*}{Proof.}
The first point is a direct consequence of the formulae contained in
Corollary \ref{cor_sigma_constant}. For the second point, the
existence of the invariant law is given by a general theorem (see
{\it e.g.} \cite{daprato}, Theorem 8.6.3 p.163) while the equality
$D_+X_t=-D_-X_t$ comes once again from the formulae contained in
Corollary \ref{cor_sigma_constant}. \qed
\end{pf*}

\subsection{Second order derivatives and characterization of gradient diffusions}

In \cite{rt} Theorem 5.4, the authors give a very nice
generalization of Kolmogorov's result \cite{kolmo} based on an
integration by part formula from Malliavin calculus. Precisely, the
drift is this time not assumed to be time homogeneous and nor the
diffusion stationary. Their characterization requires that there
exists one reversible law in the reciprocal class of the diffusion.
In our case, we are also able to characterize a larger class of
Brownian diffusions. However this further needs to use second order
stochastic derivatives. The main result of our paper is the
following theorem:
\begin{thm}\label{main}
Let $X$ be given by (\ref{constant+hom}), verifying assumption (H), such
that $b\in C^{2}(\R^d)$ with bounded
derivatives, and such that for all $t\in(0,T)$ the second
order derivatives of $\nabla \log p_t$ are bounded. We then have the
following equivalence:
\begin{equation}\label{}
D_+^2X_t=D_-^2X_t \ \mbox{  for almost all }t\in(0,T) \quad
\Longleftrightarrow \quad \text{$b$ is a gradient}.
\end{equation}
\end{thm}

\begin{rem}
{\rm
\begin{enumerate}
\item
Saying that $b$ is a gradient means that we can write $b=\nabla U$
for a certain potential $U:\R^d\to\R$. It
is equivalent, by Poincar\'e lemma, to verify that $G=\partial_x
b-(\partial_x b)^*$ is identically zero.
\item When $d=1$, that is when $X$ is a {\it one-dimensional} Brownian diffusion,
the equality $D_-^2X-D_+^2X=0$ is always verified, see Lemma
\ref{div}.
\item The proof we propose here is entirely based on probabilistic arguments.
A more "classical" strategy for proving that $G\equiv 0$ when
$D_-^2X=D_+^2X$ would use the fact that we then have ${\rm
div}(p_tG_i)=0$ for any index $i$ and any time $t\in (0,T)$ (see
Lemma \ref{div}). For instance, when $d=2$, this system of
equalities reduces to $(\partial_1 b_2-\partial_2 b_1)p_t=c$ on
$\R^2$, $c$ denoting a constant. It is then not difficult to deduce
that $\partial_1 b_2=\partial_2 b_1$. In particular, $b$ is a
gradient. On the other hand this method seems hard to adapt in
higher dimensions. In particular, it seems already difficult to
integrate ${\rm div}(p_tG)=0$ when $d=3$.
\end{enumerate}
}
\end{rem}

First of all, we need the following technical lemma which gives a
remarkable identity for $D_+^2 X-D_-^2X$:

\begin{lem}\label{div}
Let $X$ be given by (\ref{constant}), verifying assumption (H), such
that $b\in C^{1,2}([0,T]\times \R^d)$ with bounded
derivatives, and such that for all $t\in(0,T)$ the second
order derivatives of $\nabla \log p_t$ are bounded. Therefore for
any $i=1,\ldots,n$:
\begin{equation}\label{D_-^2-D_+^2i}
    (D_-^2X_t-D_+^2X_t)^i=\frac{\div(p_tG_i)}{p_t}.
\end{equation}
Recall that $G=(\partial_xb)-(\partial_xb)^*$, {\it i.e.}
$G_i^j=\partial_ib^j-\partial_jb^i$.
\end{lem}
Let us stress that the expression we obtain in (\ref{D_-^2-D_+^2i})
is the key point of our proof of Theorem \ref{main}, and that it is valid for diffusions
of the type (\ref{constant}) and not only of the type (\ref{constant+hom}).
\begin{pf*}{Proof.}
We have, by Proposition \ref{nelson-composition}:
\begin{equation}\label{}
    D_+^2X_t=D_+b(t,X_t)=\left(\partial_t b+ (\partial_xb) b
    +\frac{\sg^2}{2}\Delta b\right)(t,X_t),
\end{equation}
and
\begin{eqnarray*}
  D_-^2X_t &=& D_-\left(b-\sg^2 \frac{\nabla p_t}{p_t}\right)(t,X_t) \\
   &=& \left[\partial_t b+ (\partial_xb) b
    -\frac{\sg^2}{2}\Delta b-\sg^2 \partial_t\frac{\nabla p_t}{p_t}-\sg^2
    (\partial_xb)
\frac{\nabla p_t}{p_t} \right.\\
   & &\left. -\sg^2 \left(\partial_x\frac{\nabla
   p_t}{p_t}\right) b +\sg^4\left(\partial_x\frac{\nabla p_t}{p_t}\right) \frac{\nabla
   p_t}{p_t}+\frac{\sg^4}{2} \Delta \frac{\nabla p_t}{p_t}\right](t,X_t).
\end{eqnarray*}
With the Fokker-Planck equation
$\partial_tp_t=-\div(p_tb)+\frac{\sg^2}{2}\Delta p_t$ in mind, we
can write:
\begin{equation}\label{}
\partial_t\frac{\nabla p_t}{p_t}=\nabla\frac{\partial_t p_t}{p_t}=\nabla \left(-\div b
+\frac{\langle-b,\nabla p_t\rangle+\frac{\sg^2}{2}\Delta
p_t}{p_t}\right).
\end{equation}
Therefore: $$D_-^2X_t-D_+^2X_t = (\sg^2 A+\sg^4 B)(t,X_t)$$ with
\begin{eqnarray*}
A &=& -\Delta b+\nabla \div b- (\partial_xb) \frac{\nabla p_t}{p_t}
+ \nabla \frac{\langle b,\nabla p_t\rangle}{p_t}-\left(
\partial_x\frac{\nabla
   p_t}{p_t}\right) b, \\
 B  &=&
\left(\partial_x\frac{\nabla p_t}{p_t}\right) \frac{\nabla
   p_t}{p_t}
+\frac{1}{2} \,\Delta \frac{\nabla p_t}{p_t}-\frac{1}{2} \, \nabla
\frac{\Delta p_t}{p_t}.
\end{eqnarray*}

Let us simplify $A$. By the Leibniz rule we have:
$$
\nabla \frac{\langle b,\nabla p_t\rangle}{p_t} = (\partial_xb)^*
\frac{\nabla p_t}{p_t} +\left(\partial_x\frac{\nabla
   p_t}{p_t}\right)^* b.
$$
Since $p_t\in C^2$, the Schwarz lemma yields
$\left(\partial_x\frac{\nabla
   p_t}{p_t}\right)^*=\left(\partial_x\frac{\nabla
   p_t}{p_t}\right)$. Thus
$$A=-\Delta b+\nabla \div b+ G\
\frac{\nabla p_t}{p_t},$$ from which we deduce $$
A^i=\frac{\div(p_tG_i)}{p_t}.$$

Let us simplify $B$. We have:
$$
2\left[\left(\partial_x\frac{\nabla p_t}{p_t}\right) \frac{\nabla
   p_t}{p_t}\right]^i=2\sum_j\partial_i\left(\frac{\partial_j p_t}{p_t}\right) \frac{\partial_j p_t}{p_t}=\partial_i\sum_j\left(\frac{\partial_j
   p_t}{p_t}\right)^2.
$$
But, again par the Schwarz lemma:
$$
\left[\Delta \frac{\nabla p_t}{p_t}\right]^i=\sum_j\partial_j^2
\frac{\partial_i p_t}{p_t}= \partial_i\sum_j
\partial_j \left(\frac{\partial_j
   p_t}{p_t}\right).
$$
We then deduce that $B=0$, which concludes the proof. \qed

\end{pf*}

Now, we go back to the proof of Theorem \ref{main}. In order to simplify the
exposition, in the sequel we assume without loss of generality that $\sigma=1$.
Let $\gamma:\R^d\rightarrow\R^d$ be a
bounded Lipschitz function and $X^\ep$, for $\ep >0$, be the unique
solution of
\begin{equation}\label{X_e}
dX^\ep_t=(b+\ep\gamma)(X_t^\ep)dt+ dW_t,\quad t\in[0,T],\quad
X^\ep_0=X_0\in{\rm L}^2(\Omega).
\end{equation}
Before proving Theorem 6, we need the following lemma, stated and proved in \cite{FLLLT},
Proposition 3.1:
\begin{lem}\label{flllt}
Let $\phi:C[0,T]\to\R$ be a measurable function such that
$E[\phi(X)^2]$ is finite. Then the following equality holds:
\begin{equation}\label{flllt-lm}
\frac{\partial}{\partial \ep}\,
E[\phi(X^\ep)]_{|_{\ep=0}}=E\left[\phi(X)\int_0^T
\langle\gamma(X_s),dW_s\rangle\right].
\end{equation}
\end{lem}
\begin{pf*}{Proof.}
For the sake of completeness, let us briefly recall how the authors obtain
(\ref{flllt-lm}). We can write $
E[\phi(X^\ep)]=E^{\Q^\ep}[(Z^\ep)^{-1},\phi(X^\ep)] $ with
$d\Q^\ep/d\bP=Z^\ep$, where
$$
\begin{array}{lll}
Z^\ep&=&{\rm exp}\left(-\ep\int_0^T
\langle\gamma(X^\ep_s),dW_s\rangle- \frac{\ep^2}{2}\int_0^T
|\gamma(X^\ep(s))|^2ds\right)\\
&=&{\rm exp}\left(-\ep\int_0^T
\langle\gamma(X^\ep_s),dW^\ep_s\rangle+ \frac{\ep^2}{2}\int_0^T
|\gamma(X^\ep(s))|^2ds\right)
,
\end{array}
$$
and $W^\ep_t= W_t+\ep\int_0^t \gamma(X^\ep_s)ds$. Note
that, under $\Q^\ep$, $W^\ep$ is a Brownian motion by Girsanov
theorem. In particular the law of $(X^\ep,W^\ep)$ under $\Q^\ep$ is
the same as the law of $(X,W)$ under $\bP$. Consequently, $
E[\phi(X^\ep)]=E[(Z^\ep)^{-1}\phi(X)]. $ Equality (\ref{flllt-lm})
follows now easily by Lebesgue bounded convergence. \qed
\end{pf*}

Now, we go back to the proof of Theorem \ref{main}:
\begin{pf*}{Proof.}\label{}
If $b$ is a gradient, then for any $i\in\{1,\cdots,d\}$, $G_i=0$. So
Lemma \ref{div} yields $D_-^2X_t-D_+^2X_t=0$.

Conversely, assume that $D_-^2X_t-D_+^2X_t=0$ for any $t\in(0,T)$.
Let $i\in\{1,\cdots,d\}$, $\ep\geq 0$, and $X^\ep$ be the diffusion
process defined by (\ref{X_e}) with $\gamma=G_i$. We denote by $\cal
L_\ep$ the infinitesimal generator of $X^\ep$, considered as a
$(L^2(\R^d),\langle\cdot,\cdot\rangle)$ operator. For simplicity,
$\cal L =\cal L_0$ will denote the generator of $X
= X^0$. It is well-known that the adjoint $\cal L_\ep^*$ of
$\cal L_\ep$ writes
\begin{equation}\label{generator}
  \cal L^*_\ep = -\div [(b+\ep G_i)\ \cdot\ ] +\frac{1}{2}\Delta\ .
\end{equation}
Let $f\in C^{\infty}_0(\R^d)$. The Dynkin formula for $X$ reads:
\begin{equation}\label{}
    E[f(X_t)]-f(x)=E\left[\int_0^t \cal L
    f(X_s)ds\right].
\end{equation}
But
\begin{eqnarray}\label{3=}
E\left[\int_0^t \cal L
    f(X_s)ds\right] & = & \int_0^t\int_{\R^d}\cal
    L f(y)p_s(y)dyds \nonumber \\
 & = & \int_0^t\int_{\R^d}
    f(y)\cal L^*p_s(y)dy ds \nonumber\\
 & = & \int_0^t E\left[f(X_s)\frac{\cal
 L^*p_s(X_s)}{p_s(X_s)}\right]ds.
\end{eqnarray}
Since for all $s\in(0,T)$, $\frac{\div(p_s G_i)}{p_s}(X_s)=0$ a.s.,
we deduce from (\ref{3=}) and (\ref{generator}) that:
$$E\left[\int_0^t \cal L
    f(X_s)ds\right]=\int_0^t E\left[f(X_s)\frac{\cal
 L_\ep^*p_s(X_s)}{p_s(X_s)}\right]ds=E\left[\int_0^t \cal L_\ep
    f(X_s)ds\right].$$
Therefore:
\begin{equation}\label{}
    E[f(X_t)]-f(x)=E\left[\int_0^t \cal L_\ep
    f(X_s)ds\right].
\end{equation}
So the process $M$ defined by
$$M_t = f(X_t)-f(x)-\int_0^t \cal
L_\ep f(X_s)ds$$ is a $(\scr P^W,\bP)$-martingale (recall that we decided
to note ${\scr P}^W_t$ the $\sigma$-field generated by $W_s$ for
$s\in[0,t]$, see section \ref{s21}). Indeed, by the
Markov property applied to $X$, we can write
$$
E(M_t-M_s|\scr{P}^W_s)=E_{X_s}\left(f(X_{t-s})-f(x)-\int_0^{t-s}\cal
L_\ep f(X_s)ds\right)=0.
$$
Thus the law of $X$ solves the martingale problem associated with
the Markov diffusion $X^\ep$. But $b$ has linear growth and since
the second order derivatives of $b$ are bounded it is also the case
for $G_i$ and so for $b+\ep G_i$. This allows to apply the Stroock-Varadhan theorem (see
{\it e.g.} \cite[Th 24.1 p.170]{rt}) which establishes the existence
and uniqueness of solutions for the martingale problem. Therefore
$X$ and $X^\ep$ have the same law. As a consequence, for any
measurable function $\phi:C[0,T]\to\R$ such that
$E[\phi(X)^2]<\infty$, the function $\ep\mapsto E[\phi(X^\ep)]$ is
constant. We now apply Lemma \ref{flllt} with $\gamma=G_i$. So, we
have:
\begin{equation}\label{eq24}
E\left[\phi(X)\int_0^T \langle G_i(X_s),dW_s\rangle\right]=0.
\end{equation}
Let $\Q$ be the equivalent probability to $\bP$ given by Girsanov
theorem applied to $X$. In particular, $X$ is a Brownian motion
under $\Q$ and $\eta= d\Q/d\bP\in\mathscr{F}^{X}_T$. Thanks
to (\ref{eq24}), we have $E^\Q\left[\phi(X)\eta^{-1}\int_0^T\langle
G_i(X_s),dW_s\rangle \right]=0$. Thus, since $\phi$ is arbitrary,
Lemma 1.1.3. in \cite{n} shows that $\int_0^T \langle
G^i(X_s),dW_s\rangle=0$ $\Q$-a.s. and then also $\bP$-a.s. Then
$G_i(X)\equiv 0$ by It\^{o} isometry (under $\bP$) and, since
$\mathscr{L}(X_t)$ has a positive density for any $t\in(0,T)$, we
finally have $G\equiv 0$. This concludes the proof. \qed

\end{pf*}

\section{A remark on the complex stochastic embedding of the Newton equation}
\label{}

From $D_\pm$ one can construct a {\it complex} stochastic derivative
\begin{equation}\label{cald}
\cal D=\frac{D_++D_-}{2}+i\frac{D_+-D_-}{2}
\end{equation}
which extends on stochastic processes the classical derivative
operator $\frac{d}{dt}$. Moreover, and contrary to $D_+$ and $D_-$,
the operator $\cal{D}$ has the following natural but however
remarkable property:
\begin{pro}
For $X$ given by (\ref{X}) and verifying assumption (H), we have:
\begin{equation}\label{}
    \cal{D} X_t=0\ \mbox{ for any $t\in(0,T)$} \Longleftrightarrow X \text{ is a
 constant process on $[0,T]$.}
\end{equation}
\end{pro}
\begin{pf*}{Proof.}
The condition $\cal{D}X_t=0$ is equivalent to $D_+X_t=D_-X_t=0$.
Thus the forward drift and the backward drift are zero. So $X$ is a
$(\scr P^X_t)$ and $(\scr F^X_t)$-martingale. We can then use the
arguments of Nelson in the proof of Theorem 11.11 in \cite{nelson}
which allow to conclude. \qed
\end{pf*}
Using (\ref{cald}) and extending $\cal D$ by $\C$-linearity to
complex Brownian diffusions, we can easily compute
\begin{equation}\label{4bis}
\cal D^2=\frac{D_+D_-+D_-D_+}{2}+i\frac{D^2_+-D^2_-}{2}.
\end{equation}
Let us remark that the real part of $\cal D^2$ coincides with the
notion of {\it mean acceleration} introduced by Nelson
\cite{nelson}, equality (11.15), for which he had conjectured that it is the
more relevant quantity describing an acceleration on Brownian
diffusions.

Now, as an example, let us consider the analog of the Newton
equation
\begin{equation}\label{newton-det}
\frac{d^2x}{dt^2}=-\nabla U(x)
\end{equation}
using this new derivative $\cal D$ acting on Brownian diffusions.
More precisely, assume that the Brownian diffusion $X$ given by
(\ref{constant+hom}) verifies, for any $t\in[0,T]$:
\begin{equation}\label{newton}
\cal D^2 X_t=-\nabla U(X_t).
\end{equation}
Equation (\ref{newton}) is called the {\it stochastic embedded
equation} of the Newton equation with respect to the extension $\cal
D$ (see \cite{cd}). This embedded equation contains the
deterministic ordinary differential equation, as an equation written
in the sense of distributions of the Schwartz theory is an extension
of the initial ordinary or partial differential equation.

As we said, (\ref{newton}) admits at least $t\mapsto x_t$ verifying
(\ref{newton-det}) as solution. But, what can we say about
uniqueness? If not, what can we say about the other solutions?

First, if $X$ satisfies (\ref{newton}), we must have
$D^2_+X_t=D^2_-X_t$ for any $t\in[0,T]$, see (\ref{4bis}). Moreover,
it is proved in \cite{cd} that, under some regularity conditions, if
one searches solutions of (\ref{newton}) {\it in the class of
gradient diffusions} of the type (\ref{constant+hom}), the density $p_t$
of the solution $X_t$ is characterized via the Schr\"{o}dinger
equation. Our Theorem \ref{main}, which was conjectured at the end
of the note \cite{cd}, shows that the solutions of (\ref{newton})
are forced to have a gradient function as drift coefficient.

\end{document}